\documentclass[12pt]{amsart}
\input{epsf}
\usepackage{amssymb,latexsym,cite}
\topskip  \numberwithin{equation}{section}
\newtheorem{theorem}{Theorem}[section]
\newtheorem{proposition}[theorem]{Proposition}

\newtheorem{remark}[theorem]{Remark}

\begin{document}

\pagenumbering{arabic}
\pagestyle{headings}

\newcommand{\DPB}[4]{P\beta_{#1}^{(#2)}(#3,#4)}

\title{Representations by probabilistic Frobenius-Euler and degenerate Frobenius-Euler polynomials}
\author{Taekyun Kim*}
\address{Department of Mathematics, Kwangwoon University, Seoul 139-701, Republic of Korea}
\email{tkkim@kw.ac.kr}

\author{Dae San Kim}
\address{Department of Mathematics, Sogang University, Seoul 121-742, Republic of Korea}
\email{dskim@sogang.ac.kr}

\thanks{ * corresponding author}

\begin{abstract}
Let $Y$ be a random variable whose moment generating function exists in a neighborhood of the origin. The aim of this paper is to represent arbitrary polynomials in terms of probabilistic Frobenius-Euler polynomials associated with $Y$ and probabilistic degenerate Frobenius-Euler polynomials associated with $Y$, and more generally of their higher-order counterparts. We derive explicit formulas with the help of umbral calculus and illustrate our results in the case of several random variables $Y$.
\end{abstract}
\maketitle

\markboth{\centerline{\scriptsize Representations by probabilistic Frobenius-Euler and degenerate Frobenius-Euler polynomials}}
{\centerline{\scriptsize Taekyun Kim and Dae San Kim}}

\section{Introduction and preliminaries}

Let $Y$ be a random variable whose moment generating function exists in a neighborhood of the origin (see \eqref{1a}). The study of degenerate versions of special polynomials and numbers, originating with Carlitz's work on degenerate Bernoulli and Euler polynomials \cite{4}, has seen a resurgence of interest recently \cite{6,16,18,21,22}. These degenerate versions are not only limited to special numbers and polynomials but also extended to transcendental functions and umbral calculus \cite{15,19,20}. Similarly, probabilistic extensions of special polynomials and numbers have been extensively researched \cite{1,2,3,6,16,21,22,23,32}. \par
This paper investigates the problem of representing arbitrary polynomials in terms of probabilistic Frobenius-Euler polynomials, $H_{n}^{Y}(x|u)$ (see \eqref{23a}), and probabilistic degenerate Frobenius-Euler polynomials, $h_{n,\lambda}^{Y}(x|u)$ (see \eqref{25a}) with the help of umbral calculus (see Theorems 3.1 and 3.3). We also address the problem of expressing arbitrary polynomials in terms of their higher-order counterparts $H_{n}^{Y,(r)}(x|u)$ and $h_{n,\lambda}^{Y,(r)}(x|u)$  (see Theorems 4.1 and 4.2). The contribution of this paper is the derivation of such formulas which have potential applications to finding many interesting polynomial identities. Some of the previous works related to our results are \cite{5,11,12,13,14,17,24,26}. As examples, we express $(x)_{n}$ and $x^{n}$ as linear combinations of $H_{k}^{Y}(x|u)$ and $h_{k,\lambda}^{Y}(x|u)$, for some discrete and continuous random variables $Y$ (see Section 5). This requires explicit expressions for the probabilistic Stirling numbers of the first kind associated with $Y$, $S_{1}^{Y}(n,k)$, and the probabilistic degenerate Stirling numbers of the first kind associated with $Y$, $S_{1,\lambda}^{Y}(n,k)$ \cite{16}. Crucially, $S_{1}^{Y}(n,k)$ and $S_{2}^{Y}(n,k)$, and $S_{1,\lambda}^{Y}(n,k)$ and $S_{2,\lambda}^{Y}(n,k)$, satisfy orthogonality and inverse relations (Propositions 1.1 and 1.2), which are, as inversions are needed, essential for our problems. In contrast, the definitions of $S_{1}^{Y}(n,k)$ in \cite{2} and $S_{1,\lambda}^{Y}(n,k)$ in \cite{21}, based on cumulant generating functions, respectively together with $S_{2}^{Y}(n,k)$ and $S_{2,\lambda}^{Y}(n,k)$, do not possess these properties. The examples in Section 5 rely on explicit computations of $S_{1}^{Y}(n,k)$ and $S_{1,\lambda}^{Y}(n,k)$ in \cite{16}, for several discrete and continuous random variables $Y$. Consider, for example, the Poisson random variable $Y$ with parameter $\alpha >0$, whose probability mass function is given by (see \cite{30})
\begin{equation*}
p(i)=e^{-\alpha}\frac{\alpha^{i}}{i!}, \quad i=0,1,2,\dots.
\end{equation*}
Then our results in this case for $(x)_{n}$ are as follows:
\begin{align*}
&(x)_{n}=\sum_{r=0}^{n}\Big\{\sum_{l=r}^{n}\frac{1}{\alpha^{l}}S_{1}(l,r)S_{1}(n,l)+\frac{n}{1-u}\sum_{l=r}^{n-1}\frac{1}{\alpha^{l}}S_{1}(l,r)S_{1}(n-1,l)\Big\}H_{r}^{Y}(x|u), \\
&(x)_{n}=\sum_{r=0}^{n}\Big\{\sum_{l=r}^{n}\frac{1}{\alpha^{l}}S_{1,\lambda}(l,r)S_{1}(n,l)+\frac{n}{1-u}\sum_{l=r}^{n-1}\frac{1}{\alpha^{l}}S_{1,\lambda}(l,r)S_{1}(n-1,l)\Big\}h_{r,\lambda}^{Y}(x|u).
\end{align*}
Let $B_{n}(x)$ be Bernoulli polynomials given by $\frac{t}{e^{t}-1}e^{xt}=\sum_{n=0}^{\infty}B_{n}(x)\frac{t^{n}}{n!}$. Then, for any polynomial $p(x)$ of degree $n$, we have the following formula:
\begin{equation}
p(x)=\sum_{k=0}^{n}a_{k}B_{k}(x),\quad a_{k}=\frac{1}{k!}\int_{0}^{1}p^{(k)}(x) dx,\quad(k=0,1,\dots,n), \label{0a}
\end{equation}
where $p^{(k)}(x)=(\frac{d}{dx})^{k}p(x)$. \\
In \cite{17-1}, applying \eqref{0a} to $p(x)=\sum_{k=1}^{n-1}\frac{1}{k(n-k)}B_{k}(x)B_{n-k}(x)$, we obtained the following identity
\begin{equation}
\sum_{k=1}^{n-1}\frac{B_{k}(x)B_{n-k}(x)}{k(n-k)}=\frac{2}{n}\sum_{k=0}^{n-2}\frac{1}{n-k}\binom{n}{k}B_{n-k}B_{k}(x)+\frac{2}{n}H_{n-1}B_{n}(x), \label{-1a}
\end{equation}
where $n \ge 2$, and $H_{n}=1+\frac{1}{2}+ \cdots +\frac{1}{n}$. Substituting
$x=0$ and $x=\frac{1}{2}$ into equation \eqref{-1a} yields, respectively, the Miki's identity (see \cite{25}) and the Faber-Pandharipande-Zagier (FPZ) identity (see \cite{9}).
In contrast to the considerably more complex existing proofs, our derivation of these identities relies on the remarkably straightforward formula \eqref{0a}, utilizing only derivatives and integrals of the given polynomials. For Miki's identity, Gessel \cite{10} used two distinct expressions for Stirling numbers of the second kind $S_{2}\left(n,k\right)$, Shiratani and Yokoyama \cite{31} adopted $p$-adic analysis, and Miki \cite{25} employed a formula for the Fermat quotient $\frac{a^{p}-a}{p}$ modulo $p^{2}$.
Similarly, Dunne and Schubert \cite{8} derived the FPZ identity using asymptotic expansions of special polynomials arising from quantum field theory. Zagier also provided a proof in the appendix of \cite{9}. It's worth noting that Faber and Pandharipande initially conjectured relations between Hodge integrals in Gromov-Witten theory, which necessitated the FPZ identity, in 1998. \par
The outline of this paper is as follows. In Section 1, we recall some necessary facts that are needed throughout this paper. In Section 2, we go over umbral calculus briefly. In Section 3, we derive formulas expressing arbitrary polynomials in terms of the probabilistic Frobenius-Euler polynomials associated with $Y$, $H_{n}^{Y}(x|u)$, and the probabilistic degenerate Frobenius-Euler polynomials associated with $Y$, $h_{n,\lambda}^{Y}(x|u)$. In Section 4, we derive formulas representing arbitrary polynomials in terms of their higher-order counterparts, namely $H_{n}^{Y,(r)}(x|u)$ and $h_{n,\lambda}^{Y,(r)}(x|u)$. In Section 5, we illustrate our results when $Y$ are Bernoulli, Poisson, geometric and exponential random variables. Finally, we conclude our paper in Section 6. As general references of this paper, the reader may refer to \cite{7,27,28,29,30}.

Let $Y$ be a random variable whose moment generating function exists in a neighborhood of the origin:
\begin{equation}
E[e^{Yt}]=\sum_{n=0}^{\infty}E[Y^{n}]\frac{t^{n}}{n!}\,\,\, \mathrm{exists,\,\, for}\,\,|x|<r, \label{1a}
\end{equation}
for some positive real number $r$. \\
Let $(Y_{j})_{j \ge 1}$ be a sequence of mutually independent copies of the random variable $Y$, and let
\begin{equation}
S_{k}=Y_{1}+Y_{2}+\cdots +Y_{k}, \quad (k \ge 1), \quad  S_{0}=0. \label{2a}
\end{equation}
The probabilistic Stirling numbers of the second kind associated with $Y$, $S_{2}^{Y}(n,k)$, are given by (see \eqref{2a})
\begin{align}
&\frac{1}{k!}(E[e^{Yt}]-1)^{k}=\sum_{n=k}^{\infty}S_{2}^{Y}(n,k)\frac{t^{n}}{n!}, \label{3a}\\
&S_{2}^{Y}(n,k)=\frac{1}{k!}\sum_{j=0}^{k}\binom{k}{j}(-1)^{k-j}E[S_{j}^{n}]. \nonumber
\end{align}
From the definition in \eqref{3a}, it is immediate to see that
\begin{equation}
S_{2}^{Y}(k,k)=E[Y]^{k}. \label{4a}
\end{equation} \par
Assume from now on that
\begin{equation}
E[Y] \ne 0. \label{5a}
\end{equation}
We introduce the notation:
\begin{equation}
e_{Y}(t)=E[e^{Yt}]-1. \label{6a}
\end{equation}
Then we have
\begin{equation}
\frac{1}{k!}\big(e_{Y}(t)\big)^{k}=\sum_{n=k}^{\infty}S_{2}^{Y}(n,k)\frac{t^{n}}{n!}. \label{7a}
\end{equation}
If $f(t)=\sum_{n=0}^{\infty}a_{n}\frac{t^{n}}{n!}$ is a delta series, namely $a_{0}=0$ and $a_{1} \ne 0$, then the compositional inverse $\bar{f}(t)$ of $f(t)$ satisfying $f(\bar{f}(t))=\bar{f}(f(t))=t$ exists. Note that, as $e_{Y}(t)=E[Y]t+\sum_{m=2}^{\infty}E[Y^{m}]\frac{t^{m}}{m!}$ and $E[Y] \ne 0$ (see \eqref{5a}, \eqref{6a}), $e_{Y}(t)$ is a delta series. \par
Now, we define the  probabilistic Stirling numbers of the first kind associated with $Y$ by: for $k \ge 0$,
\begin{equation}
\frac{1}{k!}\big(\bar{e}_{Y}(t)\big)^{k}=\sum_{n=k}^{\infty}S_{1}^{Y}(n,k)\frac{t^{n}}{n!}, \label{8a}
\end{equation}
where $\bar{e}_{Y}(t)$ is the compositional inverse of $e_{Y}(t)$. \\
In addition, as usual, we agree that
\begin{equation}
S_{2}^{Y}(n,k)=S_{1}^{Y}(n,k)=0,\,\, \mathrm{if}\,\, k>n\,\,\mathrm{or}\,\, k <0. \label{9a}
\end{equation}
Note that $S_{2}^{Y}(n,k)=S_{2}(n,k),\,\, S_{1}^{Y}(n,k)=S_{1}(n,k)$, for $Y=1$. \\
Here $S_{2}(n,k)$ are the Stirling numbers of the second kind defined by
\begin{align}
&x^{n}=\sum_{k=0}^{n}S_{2}(n,k)(x)_{k},\label{10a} \\
&\frac{1}{k!}(e^{t}-1)^{k}=\sum_{n=k}^{\infty}S_{2}(n,k)\frac{t^{n}}{n!}, \nonumber
\end{align}
and $S_{1}(n,k)$ are the Stirling numbers of the first kind defined as
\begin{align}
&(x)_{n}=\sum_{k=0}^{n}S_{1}(n,k)x^{k}, \label{11a} \\
&\frac{1}{k!}\big(\log(1+t)\big)^{k}=\sum_{n=k}^{\infty}S_{1}(n,k)\frac{t^{n}}{n!}, \nonumber
\end{align}
where $(x)_{n}$ are the falling factorials given by
\begin{equation}
(x)_{0}=1, \quad (x)_{n}=x(x-1)\cdots(x-n+1),\quad (n \ge 1). \label{12a}
\end{equation} \par
Using the definitions in \eqref{7a} and \eqref{8a}, one shows that $S_{2}^{Y}(n,k)$ and $S_{1}^{Y}(n,k)$ satisfy the orthogonality relations in (a) of Proposition 1.1, from which the inverse relations in (b) and (c) follow.
\begin{proposition}
The following orthogonality and inverse relations are valid for $S_{1}^{Y}(n,k)$ and $S_{2}^{Y}(n,k)$.
\begin{flalign*}
&(a)\,\, \,\sum_{k=l}^{n} S_{2}^{Y}(n,k)S_{1}^{Y}(k,l)=\delta_{n,l}, \quad \sum_{k=l}^{n} S_{1}^{Y}(n,k)S_{2}^{Y}(k,l)=\delta_{n,l}, \\
&(b)\,\, a_{n}=\sum_{k=0}^{n}S_{2}^{Y}(n,k) b_{k}\,\, \iff \,\, b_{n}=\sum_{k=0}^{n}S_{1}^{Y}(n,k)a_{k}, \\
&(c)\,\, a_{n}=\sum_{k=n}^{m}S_{2}^{Y}(k,n)b_{k} \,\, \iff \,\, b_{n}=\sum_{k=n}^{m}S_{1}^{Y}(k,n)a_{k}. &&
\end{flalign*}
\end{proposition}
Let $\lambda$ be any nonzero real number. Then $e_{\lambda}^{x}(t)$ are the degenerate exponentials defined by
\begin{equation}
e_{\lambda}^{x}(t)=(1+\lambda t)^{\frac{x}{\lambda}}=\sum_{n=0}^{\infty}(x)_{n,\lambda}\frac{t^{n}}{n!}, \quad e_{\lambda}(t)=e_{\lambda}^{1}(t),\quad (\mathrm{see}\ [15,21]), \label{13a}
\end{equation}
where $(x)_{n,\lambda}$ are the degenerate falling factorials given by
\begin{equation}
(x)_{0,\lambda}=1,\ (x)_{n,\lambda}=x(x-\lambda)\cdots(x-(n-1)\lambda),\ (n\ge 1).\label{14a}
\end{equation}
Note here that $\lim_{\lambda\rightarrow 0}e_{\lambda}^{x}(t)=e^{xt}$. \\
The probabilistic degenerate Stirling numbers of the second kind associated $Y$, $S_{2, \lambda}^{Y}(n,k)$, are defined by
\begin{align}
&\frac{1}{k!}\big(E[e_{\lambda}^{Y}(t)]-1\big)^{k}=\sum_{n=k}^{\infty}S_{2,\lambda}^{Y}(n,k)\frac{t^{n}}{n!}, \label{15a} \\
&S_{2,\lambda}^{Y}(n,k)=\frac{1}{k!}\sum_{j=0}^{k}\binom{k}{j}(-1)^{k-j}E[(S_{j})_{n,\lambda}]. \nonumber
\end{align} \par
To define the probabilistic degenerate Stirling numbers of the first kind associated with $Y$, we let
\begin{equation}
e_{Y,\lambda}(t)=E[e_{\lambda}^{Y}(t)]-1. \label{16a}
\end{equation}
Then we have
\begin{equation}
\frac{1}{k!}(e_{Y,\lambda}(t))^{k}=\sum_{n=k}^{\infty}S_{2,\lambda}^{Y}(n,k)\frac{t^{n}}{n!}. \label{17a}
\end{equation}
Noting that $e_{Y,\lambda}(t)=E[Y]t+\sum_{m=2}^{\infty}E[(Y)_{m,\lambda}]\frac{t^{m}}{m!}$ is a delta series (see \eqref{14a}), we define the probabilistic degenerate Stirling numbers of the first kind associated with $Y$, $S_{1,\lambda}^{Y}(n,k)$, by
\begin{equation}
\frac{1}{k!}\big(\bar{e}_{Y,\lambda}(t)\big)^{k}=\sum_{n=k}^{\infty}S_{1,\lambda}^{Y}(n,k)\frac{t^{n}}{n!}, \label{18a}
\end{equation}
where $\bar{e}_{Y,\lambda}(t)$ is the compositional inverse of $e_{Y,\lambda}(t)$. \\
Note that $S_{2,\lambda}^{Y}(n,k)=S_{2,\lambda}(n,k)$ and $S_{1,\lambda}^{Y}(n,k)=S_{1,\lambda}(n,k)$, for $Y=1$.
Here $S_{2,\lambda}(n,k)$ are the degenerate Stirling numbers of the second kind defined by
\begin{align}
&(x)_{n,\lambda}=\sum_{k=0}^{n}S_{2,\lambda}(n,k)(x)_{k}, \label{19a} \\
&\frac{1}{k!}(e_{\lambda}(t)-1)^{k}=\sum_{n=k}^{\infty}S_{2,\lambda}(n,k)\frac{t^{n}}{n!}, \nonumber
\end{align}
and $S_{1,\lambda}(n,k)$ are the degenerate Stirling numbers of the first kind given by
\begin{align}
&(x)_{n}=\sum_{k=0}^{n}S_{1,\lambda}(n,k)(x)_{k,\lambda}, \label{20a} \\
&\frac{1}{k!}\big(\log_{\lambda}(1+t)\big)^{k}=\sum_{n=k}^{\infty}S_{1,\lambda}(n,k)\frac{t^{n}}{n!}, \nonumber
\end{align}
where $\log_{\lambda}(t)$ are the degenerate logarithm defined by
\begin{equation}
\log_{\lambda}(t)=\frac{1}{\lambda}\big(t^{\lambda}-1\big). \label{21a}
\end{equation}
Note here that the degenerate exponential $e_{\lambda}(t)$ in \eqref{13a} and the degenerate logarithm $\log_{\lambda}(t)$ in \eqref{21a} are compositional inverses to each other so that
\begin{equation}
e_{\lambda}\big(\log_{\lambda}(t) \big)=\log_{\lambda}\big(e_{\lambda}(t) \big)=t. \label{22a}
\end{equation} \par
From \eqref{17a} and \eqref{18a}, one shows that $S_{2,\lambda}^{Y}(n,k)$ and $S_{1,\lambda}^{Y}(n,k)$ satisfy the orthogonality relations in (a) of Proposition 1.2, from which the inverse relations in (b) and (c) follow.
\begin{proposition}
The following orthogonality and inverse relations are valid for $S_{1,\lambda}^{Y}(n,k)$ and $S_{2,\lambda}^{Y}(n,k)$.
\begin{flalign*}
&(a)\,\, \sum_{k=l}^{n} S_{2,\lambda}^{Y}(n,k)S_{1,\lambda}^{Y}(k,l)=\delta_{n,l}, \quad \sum_{k=l}^{n} S_{1,\lambda}^{Y}(n,k)S_{2,\lambda}^{Y}(k,l)=\delta_{n,l},\\
&(b)\,\, a_{n}=\sum_{k=0}^{n}S_{2,\lambda}^{Y}(n,k) b_{k}\,\, \iff \,\, b_{n}=\sum_{k=0}^{n}S_{1,\lambda}^{Y}(n,k)a_{k}, \\
&(c)\,\, a_{n}=\sum_{k=n}^{m}S_{2,\lambda}^{Y}(k,n)b_{k} \,\, \iff \,\, b_{n}=\sum_{k=n}^{m}S_{1,\lambda}^{Y}(k,n)a_{k}. &&
\end{flalign*}
\end{proposition}
\noindent Notice that the Stirling numbers of both kinds and the degenerate Stirling numbers of both kinds satisfy orthogonality and inversion relations (see Propositions 1.1 and 1.2 with $Y=1$). \par
Throughout this paper, we assume that {\it{$u$ is any complex number not equal to 1.}}
The probabilistic Frobenius-Euler polynomials associated with $Y$, $H_{n}^{Y}(x|u)$, are defined by
\begin{equation}
\frac{1-u}{E[e^{Yt}]-u}\big(E[e^{Yt}] \big)^{x}=\sum_{n=0}^{\infty}H_{n}^{Y}(x|u)\frac{t^{n}}{n!}. \label{23a}
\end{equation}
More generally, for any nonnegative integer $r$, the probabilistic Frobenius-Euler polynomials of order $r$ associated with $Y$, $H_{n}^{Y,(r)}(x|u)$, are given by
\begin{equation}
\bigg(\frac{1-u}{E[e^{Yt}]-u}\bigg)^{r}\big(E[e^{Yt}] \big)^{x}=\sum_{n=0}^{\infty}H_{n}^{Y,(r)}(x|u)\frac{t^{n}}{n!}. \label{24a}
\end{equation}
The probabilistic degenerate Frobenius-Euler polynomials associated with $Y$, $h_{n,\lambda}^{Y}(x|u)$, are defined by
\begin{equation}
\frac{1-u}{E[e_{\lambda}^{Y}(t)]-u}\big(E[e_{\lambda}^{Y}(t)] \big)^{x}=\sum_{n=0}^{\infty}h_{n,\lambda}^{Y}(x|u)\frac{t^{n}}{n!}. \label{25a}
\end{equation}
More generally, for any nonnegative integer $r$, the probabilistic degenerate Frobenius-Euler polynomials of order $r$ associated with $Y$, $h_{n,\lambda}^{Y,(r)}(x|u)$, are given by
\begin{equation}
\bigg(\frac{1-u}{E[e_{\lambda}^{Y}(t)]-u}\bigg)^{r}\big(E[e_{\lambda}^{Y}(t)] \big)^{x}=\sum_{n=0}^{\infty}h_{n,\lambda}^{Y,(r)}(x|u)\frac{t^{n}}{n!}. \label{26a}
\end{equation} \par
We note that $h_{n,\lambda}^{Y}(x|u) \rightarrow H_{n}^{Y}(x|u)$, and $h_{n,\lambda}^{Y,(r)}(x|u) \rightarrow H_{n}^{Y,(r)}(x|u)$, as $\lambda$ tends to $0$.
The Frobenius-Euler polynomials $H_{n}(x|u)$ are defined by
\begin{equation}\label{27a}
\frac{1-u}{e^{t}-u}e^{xt}=\sum_{n=0}^{\infty}H_{n}(x|u)\frac{t^{n}}{n!}.
\end{equation}
When $x=0$, $H_{n}(u)=H_{n}(0|u)$ are called the Frobenius-Euler numbers.
The first few terms of $H_n(u)$ are given by:
\begin{equation*}
H_{0}(u)=1,\, H_{1}(u)=-\frac{1}{1-u},\, H_{2}(u)= \frac{1+u}{(1-u)^{2}},\,  H_{3}(u)=-\frac{u^2+4u+1}{(1-u)^3},\dots.
\end{equation*}
More generally, for any nonnegative integer $r$, the Frobenius-Euler polynomials $H_n^{(r)}(x|u)$ of order $r$ are defined by
\begin{equation}\label{28a}
\bigg(\frac{1-u}{e^t-u}\bigg)^{r}e^{xt}=\sum_{n=0}^{\infty}H_{n}^{(r)}(x|u)\frac{t^n}{n!}.
\end{equation}
The degenerate Frobenius-Euler polynomials $h_{n,\lambda}(x|u)$ are given by (see [14])
\begin{equation}\label{29a}
\frac{1-u}{e_{\lambda}(t)-u}e_{\lambda}^{x}(t)=\sum_{n=0}^{\infty}h_{n,\lambda}(x|u)\frac{t^n}{n!}.
\end{equation}
More generally, for any nonnegative integer $r$, the degenerate Frobenius-Euler polynomials $h_{n,\lambda}^{(r)}(x|u)$ of order $r$ are defined by (see [14])
\begin{equation}\label{30a}
\bigg(\frac{1-u}{e_{\lambda}(t)-u}\bigg)^{r}e_{\lambda}^{x}(t)=\sum_{n=0}^{\infty}h_{n,\lambda}^{(r)}(x|u)\frac{t^n}{n!}.
\end{equation}
We remark that $h_{n,\lambda}(x|u) \rightarrow H_{n}(x|u)$, and $h_{n,\lambda}^{(r)}(x|u) \rightarrow H_{n}^{(r)}(x|u)$, as $\lambda$ tends to $0$. When $Y=1$, $H_{n}^{Y}(x|u),\, H_{n}^{Y,(r)}(x|u),\,h_{n,\lambda}^{Y}(x|u)$, and $h_{n,\lambda}^{Y,(r)}(x|u)$ become respectively $H_{n}(x|u),\, H_{n}^{(r)}(x|u),\,h_{n,\lambda}(x|u)$, and $h_{n,\lambda}^{(r)}(x|u)$.
We recall some notations and facts about forward differences. Let $\alpha$ be any complex-valued function of the real variable $x$. Then, for any real number $a$, the forward difference $\Delta_{a}$ is given by
\begin{equation}\label{31a}
\Delta_{a}\alpha(x)=\alpha(x+a)-\alpha(x).
\end{equation}
If $a=1$, then we let
\begin{equation}\label{32a}
\Delta \alpha(x)=\Delta_{1}\alpha(x)=\alpha(x+1)-\alpha(x).
\end{equation}
In general, the $n$th oder forward differences are given by
\begin{equation}\label{33a}
\Delta_{a}^{n}\alpha(x)=\sum_{i=0}^{n}\binom{n}{i} (-1)^{n-i}\alpha(x+ia).
\end{equation}
For $a=1$, we have
\begin{equation}\label{34a}
\Delta^{n}\alpha(x)=\sum_{i=0}^{n}\binom{n}{i} (-1)^{n-i}\alpha(x+i).
\end{equation}
Finally, we recall that the Stirling numbers of the second kind $S_{2}(n,k)$ are given by
\begin{equation}\label{35a}
\frac{1}{k!}(e^{t}-1)^{k}=\sum_{n=k}^{\infty}S_{2}(n,k)\frac{t^{n}}{n!},\quad(k \ge 0).
\end{equation}

\section{Review of umbral calculus}
\vspace{0.5cm}
Here we will briefly go over very basic facts about umbral calculus. For more details on this, we recommend the reader to refer to \cite{29}.
Let $\mathbb{C}$ be the field of complex numbers. Then $\mathcal{F}$ denotes the algebra of formal power series in $t$ over $\mathbb{C}$, given by
\begin{equation*}
 \mathcal{F}=\bigg\{f(t)=\sum_{k=0}^{\infty}a_{k}\frac{t^{k}}{k!}~\bigg|~a_{k}\in\mathbb{C}\bigg\},
\end{equation*}
and $\mathbb{P}=\mathbb{C}[x]$ indicates the algebra of polynomials in $x$ with coefficients in $\mathbb{C}$. \par
The set of all linear functionals on $\mathbb{P}$ is a vector space as usual and denoted by $\mathbb{P}^{*}$. Let  $\langle L|p(x)\rangle$ denote the action of the linear functional $L$ on the polynomial $p(x)$. \par
For $f(t)\in\mathcal{F}$ with $\displaystyle f(t)=\sum_{k=0}^{\infty}a_{k}\frac{t^{k}}{k!}\displaystyle$, we define the linear functional on $\mathbb{P}$ by
\begin{equation}\label{1b}
\langle f(t)|x^{k}\rangle=a_{k}.
\end{equation}
From \eqref{1b}, we note that
\begin{equation*}
 \langle t^{k}|x^{n}\rangle=n!\delta_{n,k},\quad(n,k\ge 0),
\end{equation*}
where $\delta_{n,k}$ is the Kronecker's symbol. \par
Some remarkable linear functionals are as follows:
\begin{align}
&\langle e^{yt}|p(x) \rangle=p(y), \nonumber \\
&\langle e^{yt}-1|p(x) \rangle=p(y)-p(0), \label{2b} \\
& \bigg\langle \frac{e^{yt}-1}{t}\bigg |p(x) \bigg\rangle = \int_{0}^{y}p(u) du.\nonumber
\end{align}
Let
\begin{equation}\label{3b}
 f_{L}(t)=\sum_{k=0}^{\infty}\langle L|x^{k}\rangle\frac{t^{k}}{k!}.
\end{equation}
Then, by \eqref{1b} and \eqref{3b}, we get
\begin{displaymath}
    \langle f_{L}(t)|x^{n}\rangle=\langle L|x^{n}\rangle.
\end{displaymath}
That is, $f_{L}(t)=L$. Additionally, the map $L\longmapsto f_{L}(t)$ is a vector space isomorphism from $\mathbb{P}^{*}$ onto $\mathcal{F}$.\par
Transporting the multiplication in $\mathcal{F}$ to $\mathbb{P}^{*}$ via this isomorphism
gives an algebra structure on $\mathbb{P}^{*}$. This means that the product of $L,\,M \in \mathbb{P}^{*}$ is given by
\begin{equation*}
\langle LM|x^{n} \rangle=\sum_{k=0}^{n}\binom{n}{k}\langle L|x^{k} \rangle \langle M|x^{n-k} \rangle,
\end{equation*}
and the map $L\longmapsto f_{L}(t)$ is now an $\mathbb{C}$-algebra isomorphism.
$\mathcal{F}$ is called umbral algebra which is the algebra of $\mathbb{C}$-linear functionals on $\mathbb{P}$. The umbral calculus is characterized as the study of the umbral algebra. \par
For each nonnegative integer $k$, the differential operator $t^k$ on $\mathbb{P}$ is defined by
\begin{equation}\label{4b}
t^{k}x^n=\left\{\begin{array}{cc}
(n)_{k}x^{n-k}, & \textrm{if $k\le n$,}\\
0, & \textrm{if $k>n$.}
\end{array}\right.
\end{equation}
Hence $t^{k}x^{n}=D^{k}x^{n}$, with $D=\frac{d}{dx}$.
Extending \eqref{4b} linearly, any power series
\begin{displaymath}
 f(t)=\sum_{k=0}^{\infty}\frac{a_{k}}{k!}t^{k}\in\mathcal{F}
\end{displaymath}
gives the differential operator on $\mathbb{P}$ defined by
\begin{equation}\label{5b}
 f(t)x^n=\sum_{k=0}^{n}\binom{n}{k}a_{k}x^{n-k},\quad(n\ge 0).
\end{equation}
Hence $f(t)x^{n}=f(D)x^{n}$, with $D=\frac{d}{dx}$.
It should be observed that, for any formal power series $f(t)$ and any polynomial $p(x)$, we have
\begin{equation}\label{6b}
\langle f(t) | p(x) \rangle =\langle 1 | f(t)p(x) \rangle =f(t)p(x)|_{x=0}.
\end{equation}
Here we note that an element $f(t)$ of $\mathcal{F}$ is a formal power series, a linear functional and a differential  operator. Some notable differential operators are as follows:
\begin{align}
&e^{yt}p(x)=p(x+y), \nonumber\\
&(e^{yt}-1)p(x)=p(x+y)-p(x), \label{7b}\\
&\frac{e^{yt}-1}{t}p(x)=\int_{x}^{x+y}p(u) du.\nonumber
\end{align}

The order $o(f(t))$ of the power series $f(t)(\ne 0)$ is the smallest integer for which $a_{k}$ does not vanish. If $o(f(t))=0$, then $f(t)$ is called an invertible series. If $o(f(t))=1$, then $f(t)$ is called a delta series. \par
For $f(t),g(t)\in\mathcal{F}$ with $o(f(t))=1$ and $o(g(t))=0$, there exists a unique sequence $s_{n}(x)$ (deg\,$s_{n}(x)=n$) of polynomials such that
\begin{equation} \label{8b}
\big\langle g(t)f(t)^{k}|s_{n}(x)\big\rangle=n!\delta_{n,k},\quad(n,k\ge 0).
\end{equation}
The sequence $s_{n}(x)$ is said to be the Sheffer sequence for $(g(t),f(t))$, which is denoted by $s_{n}(x)\sim (g(t),f(t))$. We observe from \eqref{8b} that
\begin{equation}\label{9b}
s_{n}(x)=\frac{1}{g(t)}p_{n}(x),
\end{equation}
where $p_{n}(x)=g(t)s_{n}(x) \sim (1,f(t))$.\par
In particular, if $s_{n}(x) \sim (g(t),t)$, then $p_{n}(x)=x^n$, and hence
\begin{equation}\label{10b}
s_{n}(x)=\frac{1}{g(t)}x^n.
\end{equation}

It is well known that $s_{n}(x)\sim (g(t),f(t))$ if and only if
\begin{equation}\label{11b}
\frac{1}{g\big(\overline{f}(t)\big)}e^{x\overline{f}(t)}=\sum_{k=0}^{\infty}\frac{s_{k}(x)}{k!}t^{k},
\end{equation}
for all $x\in\mathbb{C}$, where $\overline{f}(t)$ is the compositional inverse of $f(t)$ such that $\overline{f}(f(t))=f(\overline{f}(t))=t$. \par

The following equations \eqref{12b}, \eqref{13b}, and \eqref{14b} are equivalent to the fact that  $s_{n}\left(x\right)$ is Sheffer for $\left(g\left(t\right),f\left(t\right)\right)$, for some invertible $g(t)$:
\begin{align}
f\left(t\right)s_{n}\left(x\right)&=ns_{n-1}\left(x\right),\quad\left(n\ge0\right),\label{12b}\\
s_{n}\left(x+y\right)&=\sum_{j=0}^{n}\binom{n}{j}s_{j}\left(x\right)p_{n-j}\left(y\right),\label{13b}
\end{align}
with $p_{n}\left(x\right)=g\left(t\right)s_{n}\left(x\right),$
\begin{equation}\label{14b}
s_{n}\left(x\right)=\sum_{j=0}^{n}\frac{1}{j!}\big\langle{g\left(\overline{f}\left(t\right)\right)^{-1}\overline{f}\left(t\right)^{j}}\big |{x^{n}\big\rangle}x^{j}.
\end{equation}

\section{Representations by probabilistic Frobenius-Euler and degenerate Frobenius-Euler polynomials}

Our goal here is to find formulas expressing arbitrary polynomials in terms of probabilistic Frobenius-Euler polynomials associated with $Y$, $H_{n}^{Y}(x|u)$, and probabilistic
degenerate Frobenius-Euler polynomials associated with $Y$, $h_{n,\lambda}^{Y}(x|u)$. \\
(a) Firstly, we treat the problem of representing arbitrary polynomials by the probabilistic Frobenius-Euler polynomials associated with $Y$.
From \eqref{23a} and \eqref{11b}, we first observe that
\begin{align}
&H_{n}^{Y}(x|u) \sim \big(g(t)=\frac{e^t -u}{1-u}, f(t)\big), \label{1c}\\
& (x)_{n} \sim (1, e^{t}-1),\label{2c}
\end{align}
where the compositional inverse of $f(t)$ is given by $\bar{f}(t)=\log E[e^{Yt}]$.
Here and in the sequel {\it{we denote $\frac{e^{t}-u}{1-u}$ simply by $g(t)$, with the understanding that $u$ is a fixed complex number not equal to 1.}} \\
From \eqref{12b}, \eqref{1c} and \eqref{2c}, we note that
\begin{equation}
f(t)H_{n}^{Y}(x|u)=nH_{n-1}^{Y}(x|u),\quad  (e^{t}-1)(x)_{n}=n(x)_{n-1},\label{3c}
\end{equation}
and hence, by \eqref{7b} and \eqref{3c}, we get
\begin{equation}
\Delta(x)_{n}=(e^{t}-1)(x)_{n}=n(x)_{n-1}. \label{4c}
\end{equation}
Here we need to observe that $\bar{f}(t)$ is a delta series. Indeed, one shows that
\begin{equation*}
\bar{f}(t)=E[Y]t+\sum_{n=2}^{\infty}\sum_{j=1}^{n}(-1)^{j-1}(j-1)!S_{2}^{Y}(n,j)\frac{t^{n}}{n!}.
\end{equation*} \par
We note from \eqref{7a}, \eqref{11a} and \eqref{23a} that
\begin{align}
\sum_{n=0}^{\infty}&\big(H_{n}^{Y}(x+1|u)-uH_{n}^{Y}(x|u) \big)\frac{t^{n}}{n!}=(1-u)e^{x \log E[e^{Yt}]} \label{5c}\\
&=(1-u) \sum_{j=0}^{\infty}x^{j}\frac{1}{j!}\Big(\log\big(1+(E[e^{Yt}]-1)\big)\Big)^{j} \nonumber \\
&=(1-u) \sum_{j=0}^{\infty} x^{j}\sum_{k=j}^{\infty}S_{1}(k,j)\frac{1}{k!}(E[e^{Yt}]-1)^{k}\nonumber\\
&=(1-u) \sum_{j=0}^{\infty} x^{j}\sum_{k=j}^{\infty}S_{1}(k,j)\sum_{n=k}^{\infty}S_{2}^{Y}(n,k)\frac{t^{n}}{n!} \nonumber \\
&=(1-u) \sum_{n=0}^{\infty}\sum_{k=0}^{n}\sum_{j=0}^{k}S_{1}(k,j)S_{2}^{Y}(n,k)x^{j}\frac{t^{n}}{n!}. \nonumber
\end{align}
Hence, from \eqref{5c} and \eqref{11a}, we get
\begin{align}
H_{n}^{Y}(x+1|u)-uH_{n}^{Y}(x|u) &=(1-u)\sum_{k=0}^{n}S_{2}^{Y}(n,k)\sum_{j=0}^{k}S_{1}(k,j) x^{j} \label{6c} \\
&=(1-u)\sum_{k=0}^{n}S_{2}^{Y}(n,k) (x)_{k}, \quad (n \ge 0). \nonumber
\end{align}
As $S_{2}^{Y}(n,0)=\delta_{n,0}$ (see \eqref{7a}), from \eqref{6c} we obtain
\begin{equation}
H_{n}^{Y}(1|u)-uH_{n}^{Y}(u)=(1-u)S_{2}^{Y}(n,0)=(1-u)\delta_{n,0}. \label{7c}
\end{equation}
Here $\delta_{n,0}$ is the Kronecker's delta. \par
Let $p(x) \in \mathbb{C}[x]$ be a polynomial of degree $n$, and let
\begin{equation}
p(x)=\sum_{k=0}^{n}a_{k}H_{k}^{Y}(x|u). \label{8c}
\end{equation}
Now, for a fixed $u \ne 1$, we consider the following:
\begin{equation}
a(x)=p(x+1)-up(x). \label{9c}
\end{equation}
Then, from \eqref{6c}, \eqref{8c} and \eqref{9c}, we have
\begin{align}
&a(x)=\sum_{k=0}^{n} a_{k}\big(H_{k}^{Y}(x+1|u)-uH_{k}^{Y}(x|u)\big) \label{10c}\\
&=(1-u)\sum_{k=0}^{n}a_{k}\sum_{j=0}^{k}S_{2}^{Y}(k,j)(x)_{j}. \nonumber
\end{align}
For any integer $r$ with $0 \le r \le n$, from \eqref{4c} and \eqref{10c} we obtain
\begin{align}
\Delta^{r}a(x)&=(1-u)\sum_{k=r}^{n}a_k\sum_{j=r}^{k}S_{2}^{Y}(k,j)\Delta^{r}(x)_{j} \label{11c} \\
&=(1-u)\sum_{k=r}^{n}a_k\sum_{j=r}^{k}S_{2}^{Y}(k,j)(j)_{r}(x)_{j-r}.\nonumber
\end{align}
Letting $x=0$ in \eqref{11c} gives
\begin{equation}\label{12c}
\frac{1}{(1-u) r!}\Delta^{r}a(x) \big\vert _{x=0}=\sum_{k=r}^{n}a_{k}S_{2}^{Y}(k,r),\,\,(0 \le r \le n).
\end{equation}
By invoking the inversion in Proposition 1.1 (c) to \eqref{12c} and recalling \eqref{9c}, we finally have
\begin{align}
a_{r}&=\frac{1}{1-u}\sum_{j=r}^{n}S_{1}^{Y}(j,r)\frac{1}{j!}\Delta^{j}a(x)\big \vert_{x=0} \label{13c} \\
&=\frac{1}{1-u}\sum_{j=r}^{n}S_{1}^{Y}(j,r)\frac{1}{j!}\Delta^{j}\big(p(x+1)-up(x) \big)\big \vert_{x=0} \nonumber \\
&=\frac{1}{1-u}\sum_{j=r}^{n}S_{1}^{Y}(j,r)\frac{1}{j!}\big(\Delta^{j}p(1)-u \Delta^{j}p(0) \big). \nonumber
\nonumber
\end{align}
An alternative expression of \eqref{13c} is given by (see \eqref{7b}, \eqref{35a})
\begin{align}
a_r&=\frac{1}{1-u}\sum_{j=r}^{n}S_{1}^{Y}(j,r)\frac{1}{j!}\Delta^{j}a(x)\big \vert_{x=0} \label{14c}\\
&=\frac{1}{1-u}\sum_{j=r}^{n}S_{1}^{Y}(j,r)\frac{1}{j!}(e^{t}-1)^{j}a(x)\big \vert_{x=0}\nonumber \\
&=\frac{1}{1-u}\sum_{j=r}^{n}S_{1}^{Y}(j,r)\sum_{k=j}^{n}S_{2}(k,j)\frac{1}{k!}a^{(k)}(x) \big \vert_{x=0} \nonumber \\
&=\frac{1}{1-u}\sum_{j=r}^{n}S_{1}^{Y}(j,r)\sum_{k=j}^{n}S_{2}(k,j)\frac{1}{k!}\big(p^{(k)}(1)-up^{(k)}(0) \big) \nonumber \\
&=\frac{1}{1-u}\sum_{k=r}^{n}\sum_{j=r}^{k}S_{1}^{Y}(j,r)S_{2}(k,j)\frac{1}{k!}\big(p^{(k)}(1)-up^{(k)}(0) \big), \nonumber
\end{align} \par
where $p^{(k)}(x)=(\frac{d}{dx})^{k}p(x)$.
Using \eqref{34a}, we have another alternative expression of \eqref{13c} which is given by
\begin{align}
a_{r}&=\frac{1}{1-u}\sum_{j=r}^{n}S_{1}^{Y}(j,r)\frac{1}{j!}\Delta^{j}a(x)\big \vert_{x=0} \label{15c}\\
&=\frac{1}{1-u}\sum_{j=r}^{n}S_{1}^{Y}(j,r)\frac{1}{j!}\sum_{i=0}^{j}\binom{j}{i}(-1)^{j-i}a(i) \nonumber \\
&=\frac{1}{1-u}\sum_{j=r}^{n}\sum_{i=0}^{j}(-1)^{j-i}\frac{1}{j!}\binom{j}{i}S_{1}^{Y}(j,r)\big(p(i+1)-up(i) \big). \nonumber
\end{align}
Summarizing, from \eqref{13c}-\eqref{15c}, we get the following theorem.
\begin{theorem}
Let $p(x) \in \mathbb{C}[x], with\,\, \mathrm{deg}\, p(x)=n$, and let $p(x)=\sum_{r=0}^{n}a_{r} H_{r}^{Y}(x|u)$. Then the coefficients $a_{r}$ are given by
\begin{align*}
a_{r}&=\frac{1}{1-u}\sum_{j=r}^{n}S_{1}^{Y}(j,r)\frac{1}{j!}\big(\Delta^{j}p(1)-u \Delta^{j}p(0) \big) \\
&=\frac{1}{1-u}\sum_{k=r}^{n}\sum_{j=r}^{k}S_{1}^{Y}(j,r)S_{2}(k,j)\frac{1}{k!}\big(p^{(k)}(1)-up^{(k)}(0) \big) \\
&=\frac{1}{1-u}\sum_{j=r}^{n}\sum_{i=0}^{j}(-1)^{j-i}\frac{1}{j!}\binom{j}{i}S_{1}^{Y}(j,r)\big(p(i+1)-up(i) \big).
\end{align*}
\end{theorem}

\begin{remark}
Let $p(x) \in \mathbb{C}[x], with \,\,\mathrm{deg}\, p(x)=n$. Write $p(x)=\sum_{k=0}^{n} a_{k}H_{k}(x|u)$. \\
As $\sum_{j=r}^{k}S_{1}(j,r)S_{2}(k,j)=\delta_{k,r}$, we recover from Theorem 3.1 the result in \textnormal{\cite{12,14}}. Namely, we have
\begin{equation*}
a_{k}=\frac{1}{(1-u)k!}(p^{(k)}(1)-u p^{(k)}(0)),\,\,\textnormal{for}\,\,k=0,1,\dots,n.
\end{equation*}
\end{remark}

(b) Secondly, we treat the problem of representing arbitrary polynomials by the probabilistic degenerate Frobenius-Euler polynomials associated with $Y$. \\
From \eqref{25a} and \eqref{11b}, we first observe that
\begin{equation}
h_{n,\lambda}^{Y}(x|u) \sim \big(g(t)=\frac{e^t -u}{1-u}, f(t)\big), \label{16c}\\
\end{equation}
where the compositional inverse of $f(t)$ is given by $\bar{f}(t)=E[e_{\lambda}^{Y}(t)]$.\\
Here we need to show that $\bar{f}(t)=\log E[e_{\lambda}^{Y}(t)]$ is a delta series. Indeed, one shows that
\begin{equation*}
\bar{f}(t)=E[Y]t+\sum_{n=2}^{\infty}\sum_{j=1}^{n}(-1)^{j-1}(j-1)!S_{2,\lambda}^{Y}(n,j)\frac{t^{n}}{n!}.
\end{equation*}
From \eqref{12b} and \eqref{16c}, we note that
\begin{equation}
f(t)h_{n,\lambda}^{Y}(x|u)=nh_{n-1,\lambda}^{Y}(x|u). \label{17c }\\
\end{equation}
From \eqref{25a} and proceeding just as in \eqref{5c}, we have
\begin{equation}
h_{n,\lambda}^{Y}(x+1|u)-uh_{n,\lambda}^{Y}(x|u)=(1-u)\sum_{k=0}^{n}S_{2,\lambda}^{Y}(n,k)(x)_{k}. \label{18c}
\end{equation}
As $S_{2,\lambda}^{Y}(n,0)=\delta_{n,0}$ (see \eqref{15a}), from \eqref{18c} we get
\begin{equation}
h_{n,\lambda}^{Y}(1|u)-uh_{n,\lambda}^{Y}(u)=(1-u)\delta_{n,0},\label{19c}
\end{equation}
where $h_{n,\lambda}^{Y}(u)=h_{n,\lambda}^{Y}(0|u)$, and $\delta_{n,0}$ is the Kronecker's delta. \par
Now, we assume that $p(x) \in \mathbb{C}[x]$ has degree $n$, and write $p(x)=\sum_{k=0}^{n}a_{k}h_{k,\lambda}^{Y}(x|u)$.
For a fixed $u \ne 1$, let $a(x)=p(x+1)-up(x)$. Then, from \eqref{18c}, we have
\begin{align}
a(x)&=\sum_{k=0}^{n}a_{k}\big(h_{k,\lambda}^{Y}(x+1|u)-uh_{k,\lambda}^{Y}(x|u)\big) \label{20c}\\
&=(1-u)\sum_{k=0}^{n}a_{k}\sum_{j=0}^{k}S_{2,\lambda}^{Y}(k,j)(x)_{j}.\nonumber
\end{align}
For any integer with $0 \le r \le n$, from \eqref{20c} and \eqref{4c} we obtain
\begin{align}
\Delta^{r}a(x)&=(1-u)\sum_{k=r}^{n}a_k \sum_{j=r}^{k}S_{2,\lambda}^{Y}(k,j)\Delta^{r}(x)_{j}\label{21c} \\
&=(1-u)\sum_{k=r}^{n}a_k \sum_{j=r}^{k}S_{2,\lambda}^{Y}(k,j)(j)_{r}(x)_{j-r}.\nonumber
\end{align}
Letting $x=0$ in \eqref{21c} gives
\begin{equation}\label{22c}
\frac{1}{(1-u) r!}\Delta^{r}a(x) \big\vert _{x=0}=\sum_{k=r}^{n}a_{k}S_{2,\lambda}^{Y}(k,r),\,\,(0 \le r \le n).
\end{equation}
By invoking the inversion in Proposition 1.2 (c) to \eqref{22c}, we obtain
\begin{align}
a_{r}&=\sum_{j=r}^{n}S_{1,\lambda}^{Y}(j,r)\frac{1}{(1-u)j!}\Delta^{j}a(x) \big \vert_{x=0} \label{23c}\\
&=\frac{1}{1-u}\sum_{j=r}^{n}S_{1,\lambda}^{Y}(j,r)\frac{1}{j!}\big(\Delta^{j}p(1)-u \Delta^{j}p(0) \big). \nonumber
\end{align}
Using \eqref{34a} and proceeding just as in \eqref{15c}, we obtain an alternative expression of \eqref{23c} which is given by
\begin{equation}
a_{r}=\frac{1}{1-u}\sum_{j=r}^{n}\sum_{i=0}^{j}(-1)^{j-i}\frac{1}{j!}\binom{j}{i}S_{1,\lambda}^{Y}(j,r)\big(p(i+1)-up(i) \big). \label{24c}
\end{equation}
From \eqref{35a} and \eqref{7b}, we get yet another expression of \eqref{23c} as follows:
\begin{equation}
a_{r}=\frac{1}{1-u}\sum_{k=r}^{n}\sum_{j=r}^{k}S_{1,\lambda}^{Y}(j,r)S_{2}(k,j)\frac{1}{k!}\big(p^{(k)}(1)-up^{(k)}(0) \big), \label{25c}
\end{equation}
where $p^{(l)}(x)=(\frac{d}{dx})^{l}p(x)$. \par
Finally, from \eqref{23c}-\eqref{25c}, we get the following theorem.
\begin{theorem}
Let $p(x) \in \mathbb{C}[x], with\,\, \mathrm{deg}\, p(x)=n$, and let $p(x)=\sum_{k=0}^{n}a_{r} h_{r,\lambda}^{Y}(x|u)$. Then the coefficients $a_{r}$ are given by
\begin{align*}
a_{r}&=\frac{1}{1-u}\sum_{j=r}^{n}S_{1,\lambda}^{Y}(j,r)\frac{1}{j!}\big(\Delta^{j}p(1)-u \Delta^{j}p(0) \big) \\
&=\frac{1}{1-u}\sum_{k=r}^{n}\sum_{j=r}^{k}S_{1,\lambda}^{Y}(j,r)S_{2}(k,j)\frac{1}{k!}\big(p^{(k)}(1)-up^{(k)}(0) \big) \\
&=\frac{1}{1-u}\sum_{j=r}^{n}\sum_{i=0}^{j}(-1)^{j-i}\frac{1}{j!}\binom{j}{i}S_{1,\lambda}^{Y}(j,r)\big(p(i+1)-up(i) \big).
\end{align*}
\end{theorem}
\begin{remark}
As $(x)_{n,\lambda} \sim \big(1, \frac{1}{\lambda}(e^{\lambda t}-1)\big)$ (see \eqref{13a}), we have
\begin{equation}
\frac{1}{\lambda}(e^{\lambda t}-1)(x)_{n,\lambda}=n(x)_{n-1,\lambda}=\frac{1}{\lambda}\Delta_{\lambda}(x)_{n,\lambda}. \label{26c}
\end{equation}
When $Y=1$, from \eqref{20c} we have
\begin{equation}
a(x)=(1-u)\sum_{k=0}^{n}a_{k}\sum_{j=0}^{k}S_{2,\lambda}(k,j)(x)_{j}=(1-u)\sum_{k=0}^{n}a_{k}(x)_{k,\lambda}. \label{27c}
\end{equation}
So, when $Y=1$, we may apply $(\frac{1}{\lambda}\Delta_{\lambda})^{r}=\frac{1}{\lambda^{r}}\Delta_{\lambda}^{r}$ to \eqref{27c}, and get
\begin{equation*}
\frac{1}{\lambda^{r}}\Delta_{\lambda}^{r}a(x)=(1-u)\sum_{k=r}^{n}a_{k}(k)_{r}(x)_{k-r,\lambda}.
\end{equation*}
Then we can proceed just as before. For the details on these, one may refer to \cite{18}.
\end{remark}

\section{Representations by probabilisitc higher-order Frobenius-Euler and degenerate Frobenius-Euler polynomials}

Our goal here is to deduce formulas expressing arbitrary polynomials in terms of probabilistic Frobenius-Euler polynomials of order $r$ associated with $Y$, $H_{n}^{Y,(r)}(x|u)$, and probabilistic degenerate Frobenius-Euler polynomials of order $r$ associated with $Y$, $h_{n,\lambda}^{Y,(r)}(x|u)$. \\
(a) Firstly, we treat the problem of representing arbitrary polynomials by the probabilistic Frobenius-Euler polynomials of order $r$ associated with $Y$.
From \eqref{24a} and \eqref{11b}, we note that
\begin{equation*}
H_{n}^{Y,(r)}(x|u) \sim (g(t)^{r}, f(t)),
\end{equation*}
where $g(t)=\frac{e^t -u}{1-u}$, and the compositional inverse of $f(t)$ is given by $\bar{f}(t)=\log E[e^{Yt}]$.
From \eqref{12b}, we have
\begin{equation}
f(t)H_{n}^{Y,(r)}(x|u)=nH_{n-1}^{Y,(r)}(x|u), \label{1d}
\end{equation}
and from \eqref{24a}, it is immediate to see that
\begin{equation}
H_{n}^{Y,(r)}(x+1|u)-uH_{n}^{Y,(r)}(x|u)=(1-u)H_{n}^{Y,(r-1)}(x|u). \label{2d}
\end{equation}
We note that \eqref{2d} is equivalent to saying that
\begin{equation}
g(t)H_{n}^{Y,(r)}(x|u)=H_{n}^{Y,(r-1)}(x|u).\label{3d}
\end{equation}
Applying \eqref{3d} $r$-times, we get
\begin{equation}
g(t)^{r}H_{n}^{Y,(r)}(x|u)=H_{n}^{Y,(0)}(x|u).\label{4d}
\end{equation}
From \eqref{4d} and using \eqref{1d}, for $0 \le k \le n$, we have
\begin{equation}
f(t)^{k}g(t)^{r}H_{n}^{Y,(r)}(x|u)=f(t)^{k}H_{n}^{Y,(0)}(x|u)=(n)_{k}H_{n-k}^{Y,(0)}(x|u). \label{5d}
\end{equation}
Here, from \eqref{3a} and \eqref{24a}, we observe that
\begin{align}
\sum_{n=0}^{\infty}H_{n}^{Y,(0)}(x|u)\frac{t^{n}}{n!}&=\big(1+E[e^{Yt}]-1\big)^{x} \label{6d}\\
&=\sum_{k=0}^{\infty}(x)_{k}\frac{1}{k!}\big(E[e^{Yt}]-1\big)^{k} \nonumber \\
&=\sum_{k=0}^{\infty}(x)_{k}\sum_{n=k}^{\infty}S_{2}^{Y}(n,k)\frac{t^{n}}{n!} \nonumber \\
&=\sum_{n=0}^{\infty}\sum_{k=0}^{n}S_{2}^{Y}(n,k)(x)_{k} \frac{t^{n}}{n!}. \nonumber
\end{align}
Thus, from \eqref{6d} and \eqref{7c}, we obtain
\begin{align}
&H_{n}^{Y,(0)}(x|u)=\sum_{k=0}^{n}S_{2}^{Y}(n,k)(x)_{k}, \label{7d} \\
&H_{n}^{Y,(0)}(0|u)=S_{2}^{Y}(n,0)=\delta_{n,0}. \nonumber
\end{align}
Now, we assume that $p(x) \in \mathbb{C}[x]$ has degree $n$, and write $p(x)=\sum_{k=0}^{n}a_{k}H_{k}^{Y,(r)}(x|u)$.
Then we observe from \eqref{5d} and \eqref{7d} that
\begin{align}
f(t)^{k}g(t)^{r}p(x)&=\sum_{l=k}^{n}a_{l}\,f(t)^{k}g(t)^{r}H_{l}^{Y,(r)}(x|u) \label{8d}\\
&=\sum_{l=k}^{n}a_{l}(l)_{k}H_{l-k}^{Y,(0)}(x|u) \nonumber \\
&=\sum_{l=k}^{n}a_{l}(l)_{k}\sum_{j=0}^{l-k}S_{2}^{Y}(l-k,j)(x)_{j}. \nonumber
\end{align}
Evaluating \eqref{8d} at $x=0$ by using \eqref{7d}, we have
\begin{align}
f(t)^{k}g(t)^{r}p(x) \big \vert _{x=0}&=\sum_{l=k}^{n}a_{l}(l)_{k}S_{2}^{Y}(l-k,0)\label{9d}\\
&=\sum_{l=k}^{n}a_{l}(l)_{k}\delta_{l,k}=k!a_{k}.\nonumber
\end{align}
Thus, from \eqref{9d}, we obtain
\begin{equation}
a_k=\frac{1}{k!}f(t)^{k}g(t)^{r}p(x)\big \vert_{x=0}= \frac{1}{k!} \langle f(t)^{k}g(t)^{r} |p(x) \rangle . \label{10d}
\end{equation}
This also follows from the observation $\big\langle g(t)^{r}f(t)^{k} | H_{n}^{Y,(r)}(x|u) \big\rangle=n!\,\delta_{n,k}.$ \par
Since $g(t)=\frac{e^{t}-u}{1-u}$, from \eqref{35a} we have
\begin{align}
g(t)^{r}&=\frac{1}{(1-u)^{r}}\sum_{i=0}^{r} \binom{r}{i}(-u)^{r-i}e^{it} \label{11d}\\
&=\frac{1}{(1-u)^{r}}\sum_{i=0}^{r}\binom{r}{i}(1-u)^{r-i}(e^{t}-1)^{i} \nonumber \\
&=\frac{1}{(1-u)^{r}}\sum_{i=0}^{r}i!\binom{r}{i}(1-u)^{r-i}\sum_{m=i}^{\infty}S_{2}(m,i)\frac{t^{m}}{m!}. \nonumber
\end{align}
By using \eqref{7b} and \eqref{35a}, the equations \eqref{11d} give three alternative expressions of \eqref{10d} as in the following:
\begin{align}
a_{k}&=\frac{1}{(1-u)^{r}k!}\sum_{i=0}^{r}\binom{r}{i}(-u)^{r-i}f(t)^{k}p(x+i) \big \vert_{x=0} \label{12d} \\
&=\frac{1}{k!}\sum_{i=0}^{r}\binom{r}{i}(1-u)^{-i}\Delta ^{i} f(t)^{k}p(x) \big \vert_{x=0}\nonumber \\
&=\frac{1}{k!}\sum_{i=0}^{r}i!\binom{r}{i}(1-u)^{-i}\sum_{m=i}^{n}S_{2}(m,i)\frac{1}{m!}f(t)^{k}p^{(m)}(x) \big \vert_{x=0}. \nonumber
\end{align}

Summarizing the results so far, from \eqref{10d} and \eqref{12d} we obtain the following theorem.
\begin{theorem}
Let $p(x) \in \mathbb{C}[x], with\,\,\mathrm{deg}\, p(x)=n$. Let $g(t)=\frac{e^{t}-u}{1-u}, \bar{f}(t)=\log E[e^{Yt}]$. Then we have
\begin{equation*}
p(x)=\sum_{k=0}^{n}a_{k}H_{k}^{Y,(r)}(x|u),
\end{equation*}
\end{theorem}
where
\begin{align*}
a_k&=\frac{1}{k!}g(t)^{r}f(t)^{k}p(x)\big \vert_{x=0}\\
&=\frac{1}{(1-u)^{r}k!}\sum_{i=0}^{r}\binom{r}{i}(-u)^{r-i}f(t)^{k}p(x+i) \big \vert_{x=0} \\
&=\frac{1}{k!}\sum_{i=0}^{r}\binom{r}{i}(1-u)^{-i}\Delta ^{i} f(t)^{k}p(x) \big \vert_{x=0} \\
&=\frac{1}{k!}\sum_{i=0}^{r}\sum_{m=i}^{n}\frac{i!}{m!}\binom{r}{i}(1-u)^{-i}S_{2}(m,i)f(t)^{k}p^{(m)}(x) \big \vert_{x=0}.
\end{align*}

(b) Secondy, we treat the problem of representing arbitrary polynomials by the probabilistic degenerate Frobenius-Euler polynomials of order $r$ associated with $Y$.
From \eqref{26a} and \eqref{11b}, we note that
\begin{equation*}
h_{n,\lambda}^{Y,(r)}(x|u) \sim (g(t)^{r}, f(t)),
\end{equation*}
where $g(t)=\frac{e^t -u}{1-u}$, and the compositional inverse of $f(t)$ is given by $\bar{f}(t)=\log E[e_{\lambda}^{Y}(t)]$.
From \eqref{12b}, we have
\begin{equation}
f(t)h_{n,\lambda}^{Y,(r)}(x|u)=nh_{n-1,\lambda}^{Y,(r)}(x|u), \label{13d}
\end{equation}
and from \eqref{26a}, it is immediate to see that
\begin{equation}
g(t)h_{n,\lambda}^{Y,(r)}(x|u)=h_{n,\lambda}^{Y,(r-1)}(x|u).\label{14d}
\end{equation}
Applying \eqref{14d} $r$-times, we get
\begin{equation}
g(t)^{r}h_{n,\lambda}^{Y,(r)}(x|u)=h_{n,\lambda}^{Y,(0)}(x|u).\label{15d}
\end{equation}
From \eqref{15d} and using \eqref{13d}, for $0 \le k \le n$, we have
\begin{equation}
f(t)^{k}g(t)^{r}h_{n,\lambda}^{Y,(r)}(x|u)=f(t)^{k}h_{n,\lambda}^{Y,(0)}(x|u)=(n)_{k}h_{n-k,\lambda}^{Y,(0)}(x|u). \label{16d}
\end{equation}
Here, from \eqref{17a}, \eqref{26a} and noting that $S_{2,\lambda}^{Y}(n,0)=\delta_{n,0}$ (see \eqref{17a}), we obtain
\begin{align}
&h_{n,\lambda}^{Y,(0)}(x|u)=\sum_{k=0}^{n}S_{2,\lambda}^{Y}(n,k)(x)_{k}, \label{17d} \\
&h_{n,\lambda}^{Y,(0)}(0|u)=S_{2,\lambda}^{Y}(n,0)=\delta_{n,0}. \nonumber
\end{align}
Now, we assume that $p(x) \in \mathbb{C}[x]$ has degree $n$, and write $p(x)=\sum_{k=0}^{n}a_{k}h_{k,\lambda}^{Y,(r)}(x|u)$.
Then we observe from \eqref{16d} and \eqref{17d} that
\begin{align}
f(t)^{k}g(t)^{r}p(x)&=\sum_{l=k}^{n}a_{l}\,f(t)^{k}g(t)^{r}h_{l,\lambda}^{Y,(r)}(x|u) \label{18d}\\
&=\sum_{l=k}^{n}a_{l}(l)_{k}h_{l-k,\lambda}^{Y,(0)}(x|u) \nonumber \\
&=\sum_{l=k}^{n}a_{l}(l)_{k}\sum_{j=0}^{l-k}S_{2,\lambda}^{Y}(l-k,j)(x)_{j}. \nonumber
\end{align}
Evaluating \eqref{18d} at $x=0$ by using \eqref{17d}, we have
\begin{align}
f(t)^{k}g(t)^{r}p(x) \big \vert _{x=0}&=\sum_{l=k}^{n}a_{l}(l)_{k}S_{2,\lambda}^{Y}(l-k,0)\label{19d}\\
&=\sum_{l=k}^{n}a_{l}(l)_{k}\delta_{l,k}=k!a_{k}.\nonumber
\end{align}
Thus, from \eqref{19d}, we obtain
\begin{equation}
a_k=\frac{1}{k!}f(t)^{k}g(t)^{r}p(x)\big \vert_{x=0}= \frac{1}{k!} \big\langle f(t)^{k}g(t)^{r} |p(x) \big\rangle . \label{20d}
\end{equation}
This also follows from the observation $\big\langle g(t)^{r}f(t)^{k} | h_{n,\lambda}^{Y,(r)}(x|u) \big\rangle=n!\,\delta_{n,k}.$ \par

Now, from \eqref{11d}, \eqref{12d} and \eqref{20d} we obtain the following theorem.
\begin{theorem}
Let $p(x) \in \mathbb{C}[x], with\,\,\mathrm{deg}\, p(x)=n$. Let $g(t)=\frac{e^{t}-u}{1-u}, \bar{f}(t)=\log E[e_{\lambda}^{Y}(t)]$. Then we have
\begin{equation*}
p(x)=\sum_{k=0}^{n}a_{k}h_{k,\lambda}^{Y,(r)}(x|u),
\end{equation*}
\end{theorem}
where
\begin{align*}
a_k&=\frac{1}{k!}g(t)^{r}f(t)^{k}p(x)\big \vert_{x=0}\\
&=\frac{1}{(1-u)^{r}k!}\sum_{i=0}^{r}\binom{r}{i}(-u)^{r-i}f(t)^{k}p(x+i) \big \vert_{x=0} \\
&=\frac{1}{k!}\sum_{i=0}^{r}\binom{r}{i}(1-u)^{-i}\Delta ^{i} f(t)^{k}p(x) \big \vert_{x=0} \\
&=\frac{1}{k!}\sum_{i=0}^{r}\sum_{m=i}^{n}\frac{i!}{m!}\binom{r}{i}(1-u)^{-i}S_{2}(m,i)f(t)^{k}p^{(m)}(x) \big \vert_{x=0}.
\end{align*}

\section{Examples}
Here we express $(x)_{n}$ and $x^{n}$ as linear combinations of probabilistic Frobenius-Euler polynomials associated with $Y$, $H_{k}^{Y}(x|u)$, and probabilistic degenerate Frobenius-Euler polynomials associated with $Y$, $h_{k,\lambda}^{Y}(x|u)$, for several discrete and continuous random variables $Y$. We use the first formulas in Theorems 3.1 and 3.3 for $(x)_{n}$, and  the second ones in Theorems 3.1 and 3.3 for $x^{n}$. For those random variables $Y$, we need the explicit computations in \cite{16}, for probabilistic Stirling numbers of the first kind, $S_{1}^{Y}(n,k)$, and probabilistic degenerate Stirling numbers of the first kind, $S_{1,\lambda}^{Y}(n,k)$.

\vspace{0.1in}
Firstly, we let $(x)_{n}=\sum_{r=0}^{n}a_{r}H_{r}^{Y}(x|u)$. Then, from Theorem 3.1, we have
\begin{align}
a_{r}&=\frac{1}{1-u}\sum_{j=r}^{n}S_{1}^{Y}(j,r)\frac{1}{j!}\big(\Delta^{j}(x+1)_{n}-u\Delta^{j}(x)_{n} \big)\big \vert_{x=0} \label{1e} \\
&=\frac{1}{1-u}\sum_{j=r}^{n}\binom{n}{j}S_{1}^{Y}(j,r)\big((x+1)_{n-j}-u(x)_{n-j} \big) \big \vert_{x=0} \nonumber \\
&=\frac{1}{1-u}\sum_{j=r}^{n}\binom{n}{j}S_{1}^{Y}(j,r)\big((1)_{n-j}-u(0)_{n-j} \big). \nonumber
\end{align}
Before proceeding further, we note that
\begin{equation}
(1)_{k}-u(0)_{k}=\left\{\begin{array}{ccc}
1-u, & \textrm{if $k= 0$},\\
1, & \textrm{if $k=1$}, \\
0, & \textrm{if $k \ge 2$}.
\end{array}\right. \label{2e}
\end{equation}
For $r$ with $0 \le r \le n-1$, from \eqref{1e} and \eqref{2e} we have
\begin{align}
a_{r}&=\frac{1}{1-u}\binom{n}{n}S_{1}^{Y}(n,r)\big((1)_{0}-u(0)_{0} \big) \label{3e} \\
&\quad +\frac{1}{1-u}\binom{n}{n-1}S_{1}^{Y}(n-1,r)\big((1)_{1}-u(0)_{1} \big) \nonumber \\
&\quad + \frac{1}{1-u}\sum_{j=r}^{n-2}\binom{n}{j}S_{1}^{Y}(j,r)\big((1)_{n-j}-u(0)_{n-j} \big) \nonumber \\
&=S_{1}^{Y}(n,r)+\frac{n}{1-u}S_{1}^{Y}(n-1,r). \nonumber
\end{align}
As we see $a_{n}=S_{1}^{Y}(n,n)=\frac{1}{E[Y]^{n}}$ from \eqref{1e}, we get
\begin{equation}
a_{r}=S_{1}^{Y}(n,r)+\frac{n}{1-u}S_{1}^{Y}(n-1,r), \quad(0 \le r \le n),\quad (\mathrm{see}\,\, \eqref{9a}). \label{4e}
\end{equation}
Thus we have shown that
\begin{equation}
(x)_{n}=\sum_{r=0}^{n}\Big\{S_{1}^{Y}(n,r)+\frac{n}{1-u}S_{1}^{Y}(n-1,r) \Big\}H_{r}^{Y}(x|u). \label{5e}
\end{equation}
Secondly, we let $x^{n}=\sum_{r=0}^{n}a_{r}H_{r}^{Y}(x|u)$. Then, from Theorem 3.1, we have
\begin{align}
a_{r}&=\frac{1}{1-u}\sum_{k=r}^{n}\sum_{j=r}^{k}S_{1}^{Y}(j,r)S_{2}(k,j)\frac{1}{k!}\Big(\Big(\frac{d}{dx}\Big)^{k}(x+1)^{n}-u\Big(\frac{d}{dx} \Big)^{k} x^{n} \Big)\Big \vert_{x=0} \label{6e} \\
&=\frac{1}{1-u}\sum_{k=r}^{n}\sum_{j=r}^{k}\binom{n}{k}S_{1}^{Y}(j,r)S_{2}(k,j)\big((x+1)^{n-k}-u x^{n-k}\big) \big \vert_{x=0} \nonumber \\
&=\frac{1}{1-u}\sum_{k=r}^{n}\sum_{j=r}^{k}\binom{n}{k}S_{1}^{Y}(j,r)S_{2}(k,j)(1-u \delta_{n,k}). \nonumber \\
&=\sum_{j=r}^{n}S_{1}^{Y}(j,r)S_{2}(n,j)+\frac{1}{1-u}\sum_{k=r}^{n-1}\sum_{j=r}^{k}\binom{n}{k}S_{1}^{Y}(j,r)S_{2}(k,j), \quad (0 \le r \le n), \nonumber
\end{align}
where we understand that the second sum is 0 when $r=n$. \par
Thus we have found that
\begin{align}
x^{n}&=\sum_{r=0}^{n}\Big\{\sum_{j=r}^{n}S_{1}^{Y}(j,r)S_{2}(n,j) \label{7e} \\
&\quad\quad\quad +\frac{1}{1-u}\sum_{k=r}^{n-1}\sum_{j=r}^{k}\binom{n}{k}S_{1}^{Y}(j,r)S_{2}(k,j) \Big\}H_{r}^{Y}(x|u). \nonumber
\end{align}
Thirdly, we let $(x)_{n}=\sum_{r=0}^{n}a_{r}h_{r,\lambda}^{Y}(x|u)$. Then, from Theorem 3.3 and proceeding just as in \eqref{1e} and \eqref{3e}, we get
\begin{equation}
(x)_{n}=\sum_{r=0}^{n}\Big\{S_{1,\lambda}^{Y}(n,r)+\frac{n}{1-u}S_{1,\lambda}^{Y}(n-1,r) \Big\}h_{r,\lambda}^{Y}(x|u). \label{8e}
\end{equation}
Fourthly, we let $x^{n}=\sum_{r=0}^{n}a_{r}h_{r,\lambda}^{Y}(x|u)$. Then, from Theorem 3.3 and proceeding just as in \eqref{6e}, we obtain
\begin{align}
x^{n}&=\sum_{r=0}^{n}\Big\{\sum_{j=r}^{n}S_{1,\lambda}^{Y}(j,r)S_{2}(n,j) \label{9e} \\
&\quad\quad\quad +\frac{1}{1-u}\sum_{k=r}^{n-1}\sum_{j=r}^{k}\binom{n}{k}S_{1,\lambda}^{Y}(j,r)S_{2}(k,j) \Big\}h_{r,\lambda}^{Y}(x|u), \nonumber
\end{align}
where we understand that the second sum is 0 when $r=n$. \par

Summarizing our results, from \eqref{5e}, \eqref{7e}, \eqref{8e} and \eqref{9e}, we obtain the following theorem.
\begin{theorem}
We have the following expressions for $(x)_{n}$ and $x^{n}$ as linear combinations of $H_{k}^{Y}(x|u)$ and $h_{k,\lambda}^{Y}(x|u)$.
\begin{align*}
&(x)_{n}=\sum_{r=0}^{n}\Big\{S_{1}^{Y}(n,r)+\frac{n}{1-u}S_{1}^{Y}(n-1,r) \Big\}H_{r}^{Y}(x|u), \\
&x^{n}=\sum_{r=0}^{n}\Big\{\sum_{j=r}^{n}S_{1}^{Y}(j,r)S_{2}(n,j) \label{7e} \\
&\quad\quad\quad +\frac{1}{1-u}\sum_{k=r}^{n-1}\sum_{j=r}^{k}\binom{n}{k}S_{1}^{Y}(j,r)S_{2}(k,j) \Big\}H_{r}^{Y}(x|u), \\
&(x)_{n}=\sum_{r=0}^{n}\Big\{S_{1,\lambda}^{Y}(n,r)+\frac{n}{1-u}S_{1,\lambda}^{Y}(n-1,r) \Big\}h_{r,\lambda}^{Y}(x|u),\\
&x^{n}=\sum_{r=0}^{n}\Big\{\sum_{j=r}^{n}S_{1,\lambda}^{Y}(j,r)S_{2}(n,j) \\
&\quad\quad\quad +\frac{1}{1-u}\sum_{k=r}^{n-1}\sum_{j=r}^{k}\binom{n}{k}S_{1,\lambda}^{Y}(j,r)S_{2}(k,j) \Big\}h_{r,\lambda}^{Y}(x|u).
\end{align*}
\end{theorem}
(a) Let $Y$ be the Bernoulli random variable. Then the probability mass function of $Y$ is given by (see \cite{30})
\begin{equation*}
p(0)=1-p,\quad p(1)=p, \quad (0< p \le 1).
\end{equation*}
Then, from \cite{16}, we have
\begin{equation}
S_{1}^{Y}(n,k)=\frac{1}{p^{n}}S_{1}(n,k),\,\, S_{1,\lambda}^{Y}(n,k)=\frac{1}{p^{n}}S_{1,\lambda}(n,k). \label{10e}
\end{equation}
Now, from Theorem 5.1 and \eqref{10e}, we obtain
\begin{align*}
&(x)_{n}=\frac{1}{p^{n}}\sum_{r=0}^{n} \Big\{S_{1}(n,r)+\frac{np}{1-u}S_{1}(n-1,r) \Big\} H_{r}^{Y}(x|u), \\
&x^{n}=\sum_{r=0}^{n}\Big\{\sum_{j=r}^{n}\frac{1}{p^{j}}S_{1}(j,r)S_{2}(n,j)  \\
&\quad\quad\quad +\frac{1}{1-u}\sum_{k=r}^{n-1}\sum_{j=r}^{k}\binom{n}{k}\frac{1}{p^{j}}S_{1}(j,r)S_{2}(k,j) \Big\}H_{r}^{Y}(x|u), \\
&(x)_{n}=\frac{1}{p^{n}}\sum_{r=0}^{n} \Big\{S_{1,\lambda}(n,r)+\frac{np}{1-u}S_{1,\lambda}(n-1,r) \Big\} h_{r,\lambda}^{Y}(x|u), \\
&x^{n}=\sum_{r=0}^{n}\Big\{\sum_{j=r}^{n}\frac{1}{p^{j}}S_{1,\lambda}(j,r)S_{2}(n,j)  \\
&\quad\quad\quad +\frac{1}{1-u}\sum_{k=r}^{n-1}\sum_{j=r}^{k}\binom{n}{k}\frac{1}{p^{j}}S_{1,\lambda}(j,r)S_{2}(k,j) \Big\}h_{r,\lambda}^{Y}(x|u).
\end{align*}

(b) Let $Y$ be the Poisson random variable with parameter $\alpha >0$. Then the probability mass function of $Y$ is given by (see \cite{30})
\begin{equation*}
p(i)=e^{-\alpha}\frac{\alpha^{i}}{i!}, \quad i=0,1,2,\dots.
\end{equation*}
Then, from \cite{16}, we have
\begin{align}
&S_{1}^{Y}(n,k)=\sum_{l=k}^{n}\frac{1}{\alpha^{l}}S_{1}(l,k)S_{1}(n,l), \,\,
S_{1,\lambda}^{Y}(n,k)=\sum_{l=k}^{n}\frac{1}{\alpha^{l}}S_{1,\lambda}(l,k)S_{1}(n,l). \label{11e}
\end{align}
Now, from Theorem 5.1 and \eqref{11e}, we obtain
\begin{align*}
&(x)_{n}=\sum_{r=0}^{n}\Big\{\sum_{l=r}^{n}\frac{1}{\alpha^{l}}S_{1}(l,r)S_{1}(n,l)+\frac{n}{1-u}\sum_{l=r}^{n-1}\frac{1}{\alpha^{l}}S_{1}(l,r)S_{1}(n-1,l)\Big\}H_{r}^{Y}(x|u), \\
&x^{n}=\sum_{r=0}^{n}\Big\{\sum_{j=r}^{n}\sum_{l=r}^{j}\frac{1}{\alpha^{l}}S_{1}(l,r)S_{1}(j,l)S_2(n,j) \\
&\quad\quad\quad\quad\quad +\frac{1}{1-u}\sum_{k=r}^{n-1}\sum_{j=r}^{k}\sum_{l=r}^{j}\binom{n}{k}\frac{1}{\alpha^{l}}S_{1}(l,r)S_{1}(j,l)S_{2}(k,j) \Big\}H_{r}^{Y}(x|u),
\end{align*}
\begin{align*}
&(x)_{n}=\sum_{r=0}^{n}\Big\{\sum_{l=r}^{n}\frac{1}{\alpha^{l}}S_{1,\lambda}(l,r)S_{1}(n,l)+\frac{n}{1-u}\sum_{l=r}^{n-1}\frac{1}{\alpha^{l}}S_{1,\lambda}(l,r)S_{1}(n-1,l)\Big\}h_{r,\lambda}^{Y}(x|u), \\
&x^{n}=\sum_{r=0}^{n}\Big\{\sum_{j=r}^{n}\sum_{l=r}^{j}\frac{1}{\alpha^{l}}S_{1,\lambda}(l,r)S_{1}(j,l)S_2(n,j) \\
&\quad\quad\quad\quad\quad +\frac{1}{1-u}\sum_{k=r}^{n-1}\sum_{j=r}^{k}\sum_{l=r}^{j}\binom{n}{k}\frac{1}{\alpha^{l}}S_{1,\lambda}(l,r)S_{1}(j,l)S_{2}(k,j) \Big\}h_{r,\lambda}^{Y}(x|u).
\end{align*}

(c) Let $Y$ be the geometric random variable with parameter $0 <p <1$. Then the probability mass function of $Y$ is given by (see \cite{30})
\begin{equation*}
p(i)=(1-p)^{i-1}p, \quad i=1,2,\dots.
\end{equation*}
Then, from \cite{16}, we have
\begin{align}
&S_{1}^{Y}(n,k)=\sum_{l=k}^{n}\binom{n}{l}(n-1)_{n-l}p^{l}(p-1)^{n-l}S_{1}(l,k), \label{12e} \\
&S_{1,\lambda}^{Y}(n,k)=\sum_{l=k}^{n}\binom{n}{l}(n-1)_{n-l}p^{l}(p-1)^{n-l}S_{1,\lambda}(l,k). \nonumber
\end{align}
Now, from Theorem 5.1 and \eqref{12e}, we obtain
\begin{align*}
&(x)_{n}=\sum_{r=0}^{n}\Big\{\sum_{l=r}^{n}\binom{n}{l}(n-1)_{n-l}p^{l}(p-1)^{n-l}S_{1}(l,r) \\
&\quad\quad\quad\quad\quad + \frac{n}{1-u}\sum_{l=r}^{n-1}\binom{n-1}{l}(n-2)_{n-1-l}p^{l}(p-1)^{n-1-l}S_{1}(l,r) \Big\}H_{r}^{Y}(x|u), \\
&x^{n}=\sum_{r=0}^{n}\Big\{\sum_{j=r}^{n}\sum_{l=r}^{j}\binom{j}{l}(j-1)_{j-l}p^{l}(p-1)^{j-l}S_{1}(l,r)S_{2}(n,j) \\
&\quad\quad\quad\quad + \frac{1}{1-u}\sum_{k=r}^{n-1}\sum_{j=r}^{k}\sum_{l=r}^{j}\binom{n}{k}\binom{j}{l}(j-1)_{j-l}p^{l}(p-1)^{j-l}S_{1}(l,r)S_{2}(k,j) \Big\}H_{r}^{Y}(x|u), \\
&(x)_{n}=\sum_{r=0}^{n}\Big\{\sum_{l=r}^{n}\binom{n}{l}(n-1)_{n-l}p^{l}(p-1)^{n-l}S_{1,\lambda}(l,r) \\
&\quad\quad\quad\quad\quad + \frac{n}{1-u}\sum_{l=r}^{n-1}\binom{n-1}{l}(n-2)_{n-1-l}p^{l}(p-1)^{n-1-l}S_{1,\lambda}(l,r) \Big\}h_{r,\lambda}^{Y}(x|u), \\
&x^{n}=\sum_{r=0}^{n}\Big\{\sum_{j=r}^{n}\sum_{l=r}^{j}\binom{j}{l}(j-1)_{j-l}p^{l}(p-1)^{j-l}S_{1,\lambda}(l,r)S_{2}(n,j) \\
&\quad\quad\quad\quad + \frac{1}{1-u}\sum_{k=r}^{n-1}\sum_{j=r}^{k}\sum_{l=r}^{j}\binom{n}{k}\binom{j}{l}(j-1)_{j-l}p^{l}(p-1)^{j-l}S_{1,\lambda}(l,r)S_{2}(k,j) \Big\}h_{r,\lambda}^{Y}(x|u).
\end{align*}

(d) Let $Y$ be the exponential random variable with parameter $\alpha > 0$. Then the probability density function of $Y$ is given by (see \cite{30})
\begin{equation*}
f(y)=\left\{\begin{array}{ccc}
\alpha e^{-\alpha y}, & \textrm{if \,\,$y \ge 0$,} \\
0, & \textrm{if\,\, $y<0$}.
\end{array}\right.
\end{equation*}
Then, from \cite{16}, we have
\begin{align}
&S_{1}^{Y}(n,k)=(-1)^{n-k}\binom{n}{k}(n-1)_{n-k}\alpha^{k},\label{13e} \\
&S_{1,\lambda}^{Y}(n,k)=\sum_{l=k}^{n}\binom{n}{l}(-1)^{n-l}(n-1)_{n-l}\alpha^{l} \lambda^{l-k}S_{2}(l,k). \nonumber
\end{align} \par
Now, from Theorem 5.1 and \eqref{13e}, we obtain
\begin{align*}
&(x)_{n}=\sum_{r=0}^{n}\Big\{(-1)^{n-r}\binom{n}{r}(n-1)_{n-r}\alpha^{r} \\
&\quad\quad\quad\quad\quad+\frac{n}{1-u}(-1)^{n-1-r}\binom{n-1}{r}(n-2)_{n-1-r}\alpha^{r} \Big\}H_{r}^{Y}(x|u), \\
&x^{n}=\sum_{r=0}^{n} \Big\{\sum_{j=r}^{n}(-1)^{j-r}\binom{j}{r}(j-1)_{j-r}\alpha^{r}S_{2}(n,j) \\
&\quad\quad\quad\quad+\frac{1}{1-u}\sum_{k=r}^{n-1}\sum_{j=r}^{k}\binom{n}{k}(-1)^{j-r}\binom{j}{r}(j-1)_{j-r}\alpha^{r}S_{2}(k,j) \Big\}H_{r}^{Y}(x|u), \\
&(x)_{n}=\sum_{r=0}^{n}\Big\{\sum_{l=r}^{n}\binom{n}{l}(-1)^{n-l}(n-1)_{n-l}\alpha^{l}\lambda^{l-r}S_{2}(l,r) \\
&\quad\quad\quad\quad \frac{n}{1-u}\sum_{l=r}^{n-1}\binom{n-1}{l}(-1)^{n-1-l}(n-2)_{n-1-l}\alpha^{l}\lambda^{l-r}S_{2}(l,r) \Big\}h_{r,\lambda}^{Y}(x|u), \\
&x^{n}=\sum_{r=0}^{n}\Big\{\sum_{j=r}^{n}\sum_{l=r}^{j}\binom{j}{l}(-1)^{j-l}(j-1)_{j-l}\alpha^{l}\lambda^{l-r}S_{2}(l,r)S_{2}(n,j) \\
&\quad\quad\quad+\frac{1}{1-u}\sum_{k=r}^{n-1}\sum_{j=r}^{k}\sum_{l=r}^{j}\binom{n}{k}\binom{j}{l}(-1)^{j-l}(j-1)_{j-l}\alpha^{l}\lambda^{l-r}S_{2}(l,r)S_{2}(k,j) \Big\}h_{r,\lambda}^{Y}(x|u).
\end{align*}

\section{Conclusion}
This paper explored representations of arbitrary polynomials as linear combinations of probabilistic Frobenius-Euler polynomials associated with $Y$, $H_{n}^{Y}(x|u)$, and probabilistic degenerate Frobenius-Euler polynomials associated with $Y$, $h_{n,\lambda}^{Y}(x|u)$, and as linear combinations of their higher-order counterparts, $H_{n}^{Y,(r)}(x|u)$ and $h_{n,\lambda}^{Y,(r)}(x|u)$. We derived explicit coefficients for these linear combinations, expressed in terms of probabilistic Stirling numbers of the first kind associated with $Y$, $S_{1}^{Y}(n,k)$, and probabilistic degenerate Stirling numbers of the first kind associated with $Y$, $S_{1,\lambda}^{Y}(n,k)$. To demonstrate our findings, we provided concrete examples by expressing $(x)_{n}$ and $x^{n}$ as linear combinations of these polynomials for some discrete and continuous random variables $Y$, utilizing established results for $S_{1}^{Y}(n,k),\,\,and\,\, S_{1,\lambda}^{Y}(n,k)$.

\vspace{0.1in}

{\bf{Acknowledgment:}} We would like to thank Jangjeon Institute for Mathematical Sciences for the support of this research.



\vspace{0.1in}

\begin{thebibliography}{9}

\bibitem{1}
J. A. Adell, \emph{Probabilistic Stirling Numbers of the Second Kind and Applications,} J. Theor. Probab. \textbf{35} (2022), 636–652. https://doi.org/10.1007/s10959-020-01050-9
\bibitem{2}
J. A. Adell, B. B{\'e}nyi, \emph{Probabilistic Stirling numbers and applications,} Aequat. Math. \textbf{98} (2024), 1627-1646.
\bibitem{3}
J. A. Adell, A. Lekuona, \emph{A probabilistic generalization of the Stirling numbers of the second kind,} J. Number Theory \textbf{194} (2019), 225-355.
\bibitem{4}
L. Carlitz, \emph{Degenerate Stirling, Bernoulli and Eulerian numbers,} Utilitas Math. \textbf{15} (1979), 51-88.
\bibitem{5}
L. Carlitz, \emph{The product of two Eulerian polynomials,} Math. Mag. \textbf{36} (1963), no. 1, 37-41.
\bibitem{6}
L. Chen, T. Kim, D. S. Kim, H. Lee, S.-H. Lee, \emph{Probabilistic degenerate central Bell polynomials,} Math. Comput. Model. Dyn. Syst. \textbf{30} (2024), no. 1, 523-542. DOI: 10.1080/13873954.2024.2358899
\bibitem{7}
L. Comtet, \emph{Advanced combinatorics,} The art of finite and infinite expansions, Revised and enlarged edition, D. Reidel Publishing Co., Dordrecht, 1974.
\bibitem{8}
G. V. Dunne, C. Schubert, \emph{Bernoulli number identities from quantum field theory and topological string theory,} Commun. Number Theory Phys. \textbf{7} (2013), no. 2, 225-249.
\bibitem{9}
C. Faber, R. Pandharipande, \emph{Hodge integrals and Gromov-Witten theory,} Invent. Math. \textbf{139} (2000), no. 1, 173-199.
\bibitem{10}
I. M. Gessel, \emph{On Miki's identities for Bernoulli numbers,} J. Number Theory \textbf{110} (2005), no. 1, 75-82.
\bibitem{11}
Y. He, S. J. Wang, \emph{New formulae of products of the Frobenius-Euler polynomials,} Adv. Differ. Equ. 2014, 2014:261, 13 pp.
\bibitem{12}
D. S. Kim, T. Kim, \emph{Some identities of higher order Euler polynomials arising from Euler basis,} Integral Transforms Spec. Funct. \textbf{24} (2013), no. 9, 734-738.
\bibitem{13}
D. S. Kim, T. Kim, \emph{Some identities of Frobenius-Euler polynomials arising from umbral calculus,} Adv. Differ. Equ. 2012, 2012:196, 10 pp.
\bibitem{14}
D. S. Kim, T. Kim, \emph{Some new identities of Frobenius-Euler numbers and polynomials,} J. Inequal. Appl. 2012, 2012:307, 10 pp.
\bibitem{15}
D. S. Kim, T. Kim, \emph{Degenerate Sheffer sequences and $\lambda$-Sheffer sequences,} J. Math. Anal. Appl. \textbf{493} (2021), 124521.
\bibitem{16}
D. S. Kim, T. Kim, \emph{Probabilisitc Stirling and degenerate Stirling numbers,} Preprint.
\bibitem{17}
D. S. Kim, T. Kim, T. Mansour, \emph{Euler basis and the product of several Bernoulli and Euler polynomials,} Adv. Stud. Contemp. Math. (Kyungshang) \textbf{24} (2014), no. 4, 535--547.
\bibitem{17-1}
D. S. Kim, T, Kim, S.-H. Lee, Y.-H. Kim, \emph{Some identities for the product of two Bernoulli and Euler polynomials,} Adv. Differ. Equ. 2012, 2012:95, 14 pp.
\bibitem{18}
T. Kim, D. S. Kim \emph{Representation by degenerate Frobenius–Euler polynomials,} Georgian Math. J. \textbf{29} (2022), no. 5, 741-754. https://doi.org/10.1515/gmj-2022-2167
\bibitem{19}
T. Kim, D. S. Kim, \emph{Degenerate Laplace transform and degenerate gamma function,} Russ. J. Math. Phys. \textbf{24} (2017), no. 2, 241-248.
\bibitem{20}
T. Kim, D. S. Kim, \emph{Note on the degenerate gamma function,} Russ. J. Math. Phys. \textbf{27} (2020), no. 3, 352–358.
\bibitem{21}
T. Kim, D. S. Kim, \emph{Probabilistic degenerate Stirling numbers of the first kind and their applications,} Eur. J. Math. \textbf{10} (2024), Article No. 73.
\bibitem{22}
T. Kim, D. S. Kim, \emph{Probabilistic degenerate Bell polynomials associated with random variables,} Russ. J. Math. Phys. \textbf{30} (2023), no. 4, 528–542.
\bibitem{23}
T. Kim, D. S. Kim, \emph{Probabilistic Bernoulli and Euler polynomials,} Russ. J. Math.
Phys. \textbf{31} (2024), no. 1, 94-105.
\bibitem{24}
T. Kim, D. S. Kim, D. V. Dolgy, S. H. Rim, \emph{Some identities on the Euler numbers arising from Euler basis polynomials,} Ars Combin. \textbf{109} (2013), 433-446.
\bibitem{25}
H. Miki, \emph{A relation between Bernoulli numbers,} J. Number Theory \textbf{10} (1978), no. 3, 297-302.
\bibitem{26}
J. Pan, F. Yang, \emph{Some convolution identities for Frobenius-Euler polynomials,} Adv. Differ. Equ. 2017, 2017:6, 16 pp.
\bibitem{27}
J. Quaintance, H. W. Gould, \emph{Combinatorial identities for Stirling numbers, The unpublished  notes of H. W. Gould, With a foreword by George E. Andrews,} World Scientific Publishing Co. Pte. Ltd., Singapore, 2016.
\bibitem{28}
J. Riordan, \emph{Combinatorial identities,} Wiley, New York, 1968.
\bibitem{29}
S. Roman, \emph{The umbral calculus,} Pure and Applied Mathematics, 111. Academic Press, Inc. [Harcourt Brace Jovanovich, Publishers], New York, 1984.
\bibitem{30}
S. M. Ross, \emph{Introduction to probability models,} Twelfth edition, Academic Press, London, 2019.
\bibitem{31}
K. Shiratani, S. Yokoyama, \emph{An application of $p$-adic convolutions,} Mem. Fac. Sci. Kyushu Univ. Ser. A \textbf{36} (1982), no. 1, 73-83.
\bibitem{32}
R. Xu, Y. Ma, T. Kim, D. S. Kim, S. Boulaaras, \emph{Probabilistic central Bell polynomials,} Fractals, to appear.


\end{thebibliography}
\end{document}